\documentclass[12pt]{amsart}
\usepackage{amsmath, amssymb, euscript, hyperref}

\newtheorem{thm}{Theorem}[section]
\newtheorem{lem}[thm]{Lemma}
\newtheorem{prop}[thm]{Proposition}
\theoremstyle{definition}
\newtheorem{defn}[thm]{Definition}

\newtheorem{rem}[thm]{Remark}

\newcommand{\blackboard}[1]{\ensuremath{\mathbb{#1}}}

\newcommand{\Z}{\blackboard{Z}}

\address{Azer Akhmedov, Department of Mathematics,
North Dakota State University,
Fargo, ND, 58108, USA}
\email{azer.akhmedov@ndsu.edu}

\begin{document}

\begin{center} {\bf \Large On the height of subgroups of $\mathrm{Homeo}_{+}(I)$} \end{center}

\vspace{0.6cm}

\begin{center} {\bf Azer Akhmedov} \end{center}

\vspace{0.7cm}

Abstract: {\Small In [5], it is proved that a subgroup of $PL_{+}(I)$ has a finite height if and only if it is solvable. We prove the ``only if" part for any subgroup of Homeo$_{+}(I)$, and present a construction which indicates a plethora of examples of solvable groups with infinite height.} \footnote { *{\it 2010 Mathematics Subject Classification.} Primary 28D05; Secondary 37E05, 20F65. \newline
} 

\vspace{1cm}

 \section{Introduction}

 A major inspiration for this paper is a beautiful geometric characterization of solvable subgroups of $PL_{+}(I)$ obtained by Collin Bleak in his Ph.D thesis: it is proved that a group $\Gamma $ of piecewise linear homeomorphisms of the closed unit interval $I$ is solvable of degree $n$ if and only it admits a (strict) tower of cardinality $n$ but not $n+1$. Here, a tower is a collection of nested open intervals such that for every interval in the collection, an element of $\Gamma $ fixes the endpoints of the interval but does not fix any inner point. A tower is called strict if no two intervals in the collection share an end.   
 
 \medskip
 
 The described characterization has been used by C.Bleak to obtain interesting structural algebraic results on subgroups of $PL_{+}(I)$ (see [5], [6], [7]). It has also been used in [1] as a major tool to prove the so-called Girth Alternative for subgroups of $PL_{+}(I)$. 
 
 \medskip
 
 The Girth Alternative remains open for the whole group $\mathrm{Homeo}_{+}(I)$, and even for subgroups such as $\mathrm{Diff} _{+}(I)$ or $\mathrm{Diff} ^2_{+}(I)$, primarily because no analogous tools are available. On the other hand, it is plausible that the tower characterization (both directions) may hold in other subgroups of $\mathrm{Homeo}_{+}(I)$ besides $PL_{+}(I)$, especially in subgroups with increased regularity; in fact, in the last section, we will prove such a result for $\mathrm{Diff} ^{\omega }_{+}(I),$ the group of analytic diffeomorphisms of the closed interval $I$. Having such an analogues characterization may provide tools to tackle several other open problems about subgroups of $\mathrm{Homeo}_{+}(I)$.  
 
 \medskip
 
 In this paper we study the notion of height for subgroups of $\mathrm{Homeo}_{+}(I)$. We prove the tower characterization of solvable subgroups in one direction and show that the other direction fails badly in the continuous category. As a byproduct of our method we obtain interesting new results related to (but not covered by) the results of [2] and [3] (See Remarks \ref{rem:solodov} and \ref{rem:discrete} respectively).   
  
 \bigskip
 
 To state the major results of this paper we need to introduce some key notions; most of the following notions are borrowed from [5].
 
 \medskip
 
 We will call the convex hull of a point in $I$ under the action of $\Gamma \leq \mathrm{Homeo}_{+}(I)$ an {\em orbital} of $\Gamma $ if this convex hull contains more than one point. We note that the orbitals are open intervals. If $g\in \Gamma $, we will refer to an orbital of the group $\langle g \rangle $ as an orbital of $g$. If an open interval $A$ is an orbital of $g$ then the pair $(A,g)$ will be called a {\em signed orbital of} $G$. $g$ will be called the {\em signature} of the signed orbital $(A,g)$.
 
 \medskip
 
 Given a set $Y$ of signed orbitals of $G$, the symbol $S_Y$ will refer to the set of signatures
of the signed orbitals in $Y$. Similarly, the symbol $O_Y$ will refer to the set of orbitals of the
signed orbitals of $Y$. We note that the set of signed orbitals of $PL_o(I)$ is a partially ordered
set under the lexicographic order induced from the partial order on subsets of $I$ (induced by
inclusion) in the first coordinate, and the left total order of the elements of $PL_o(I)$ in the
second coordinate.
  
 \medskip
 
 A {\em tower} $T$ of $G$ is a set of signed orbitals which satisfies the following two criteria.

 \medskip
 
 1. $T$ is a chain in the partial order on the signed orbitals of $G$.
 
 \medskip
 
 2. For any $A\in O_T$, $T$ has exactly one element of the form $(A,g)$.

   \medskip

 Given a tower $T$ of $G$, if $(A,g), (B,h)\in T$ then one of $A\subseteq B$ and $B\subseteq A$ holds, with equality occurring only if $g = h$ as well. Therefore, one can visualize the tower as a stack of nested levels that are always getting wider as one goes ``up" the stack.
  
 \medskip

If $A = (a,b), B = (c,d)$ are orbitals of the group $G$ so that either $c = a$ or $d = b$, then we say that an orbital $B$ shares an end with $A$. A tower $T$ will be called {\em strict} if no two distinct orbitals in $O_T$ share an end. 
 
\medskip
 
 The cardinality of the set of $O_T$ will be called the height of $T$. The {\em height} of $G$ is the supremum of the cardinalities of the towers of the group, if this supremum is finite, and will be denoted by $h(G)$. If the supremum is not finite then we will say the group $G$ has infinite height, and write $h(G) = \infty $. The {\em strict height} of $G$ is the supremum of the cardinalities of the strict towers of the group, and will be denoted by $h_{strict}(G)$; similarly, if the supremum is not finite then we will say that $G$ has infinite strict height, and write $h_{strict}(G) = \infty $.

 \medskip
 
 The major result of [5] is the following beautiful geometric characterization of solvable subgroups of $PL_o(I)$.
 
 \medskip
 
 \begin{thm} \label{thm:bleak} A subgroup $G\leq PL_o(I)$ is non-solvable if and only if $G$ admits a tower of height $n$ for any $n\geq 1$. Moreover, if $G$ admits a tower of height $n\geq 1$ then it admits a strict tower of height $n$.
 \end{thm}

 \medskip

Thus Theorem \ref{thm:bleak} implies that for non-solvable subgroups $G\leq PL_o(I)$, we have $h(G) = h_{strict}(G) = \infty $. 
 
 \bigskip

   \begin{rem} A subgroup $\Gamma \leq \mathrm{Homeo}_{+}(I)$ with $h(\Gamma ) = 1$ is necessarily Abelian. This is because a condition $h(\Gamma ) = 1$ implies that any fixed point of any non-identity element of $\Gamma $ is a global fixed point of $\Gamma $. Thus we may assume that $\Gamma $ acts freely on $(0,1)$. Then by H\"older's Theorem $\Gamma $ is Abelian. 
   \end{rem}
   
  Now we can state the first major theorem of our paper.   
   
   \medskip
   
   \begin{thm} \label{thm:easy} Let $\Gamma \leq \mathrm{Homeo}_{+}(I)$ with $h(\Gamma ) = N$. Then $\Gamma $ is solvable of solvability degree at most $N$.
 \end{thm}
 
 \medskip
 
 \begin{rem}\label{rem:hvsh'} Let us emphasize that it is easy to find a subgroup $\Gamma \leq \mathrm{Homeo}_{+}(I)$ with finite strict height but infinite height. In fact, there exists a solvable subgroup $\Gamma $ such that $h_{strict}(\Gamma ) = 1$ but $h(\Gamma ) = \infty $. Let $\Gamma $ be a subgroup of $\mathrm{Homeo}_{+}(\mathbb{R})$ generated by the two homeomorphisms $f(x) = 2x$ and $g(x) = x+1$ (here, we identify $\mathrm{Homeo}_{+}(\mathbb{R})$ with $\mathrm{Homeo}_{+}(I)$ which can be done by identifying $I = \{0\}\cup (0,1)\cup \{1\}$ with $\{-\infty \}\cup \mathbb{R}\cup \{+\infty \}$). Then $\Gamma $ is isomorphic to a metaabelian Baumslag-Solitar group $BS(1,2)$. Every element of $\Gamma $ will have either no fixed point in $(-\infty , \infty )$ or only one fixed point. Thus $h_{strict}(\Gamma ) = 1$. On the other hand, for all natural $n$, the map $gf^n(x) = 2^nx+1$ fixes the endpoints of the interval $(-\infty , -\frac{1}{2^n-1})$, thus $h(\Gamma ) = \infty $.  
 \end{rem}
 
  \medskip
 
  In light of the above remark it is interesting to know if having finite strict height still implies solvability. By expanding on the ideas of the proof of Theorem \ref{thm:easy} we obtain the following result.    
  
  \medskip
  
  \begin{thm} \label{thm:hard} Let $\Gamma \leq \mathrm{Homeo}_{+}(I)$ with $h_{strict}(\Gamma ) = N$. Then $\Gamma $ is solvable of solvability degree at most $N+1$.
 \end{thm}
 
 \medskip
 
 It is of course natural to ask whether or not the converse of Theorem \ref{thm:hard} holds. (it holds in $PL_{+}(I)$, by the result of C.Bleak). In Section 5, we will present a counterexample to the converse. More precisely, we will prove the following

 \begin{thm} \label{thm:example} There exists a solvable group $\Gamma \leq \mathrm{Homeo}_{+}(I)$ of solvability degree two generated by two homeomorphisms such that $h_{strict}(\Gamma ) = \infty $.
\end{thm}  

Despite the above negative result (Theorem \ref{thm:example}), in subgroups with higher regularity the converse of Theorem \ref{thm:hard} seems still plausible. In Section 6, we will prove the converse under very strong regularity condition, namely, for subgroups of analytic diffeomorphisms.  
 
 \medskip
 
 \begin{thm}\label{thm:analytic} If $\Gamma \leq \mathrm{Diff} ^{\omega }_{+}(I)$ is solvable then $h_{strict}(\Gamma ) = 1$.
 \end{thm}
 
 Let us emphasize that, by the result of E.Ghys [10], a solvable subgroup of $\mathrm{Diff} ^{\omega }_{+}(\mathbb{S}^1)$ is necessarily metaabelian. Moreover, all solvable subgroups of $\mathrm{Diff} ^{\omega }_{+}(\mathbb{S}^1)$ have been classified (see [8]). Indeed, Theorem \ref{thm:analytic} can be obtained as an immediate corollary of this strong classification result of L.Burslem and A.Wilkinson. However, in the last section we present another set of ideas (to prove Theorem \ref{thm:analytic}) with interesting consequences in its own right.

  \bigskip
  
  \section {Preliminary Notions}
  
  In this section, we will introduce several notions as well as quote some well known results that will be useful in the sequel.

 \medskip
 
 An orbital of a group $G\leq \mathrm{Homeo}_{+}(I)$ is called {\em inner} if it is does not share an end with $(0,1)$. An inner orbital is called {\em maximal} if it is not contained in any other inner orbital (thus the notion of maximality is defined only for inner orbitals).
 
 \medskip

 A homeomorphism $f\in \mathrm{Homeo}_{+}(I)$ will be called {\em special} if the subgroup $\langle f \rangle $ does not possess an inner orbital. A subgroup $G\leq \mathrm{Homeo}_{+}(I)$ is called {\em special} if every element of $G$ is special. 

 \medskip

 Two homeomorphisms $f, g\in \mathrm{Homeo}_{+}(I)$ are said to form {\em a crossed pair} if for some interval $(a,b)\subseteq (0,1)$, one of homeomorphisms fixes $a$ and $b$ but no other point in $[a,b]$ while the other one sends $a$ or $b$ into $(a,b)$. It is a well known result that if $f, g$ form a crossed pair then the subgroup generated by $f$ and $g$ contains a free semigroup on two generators (see, for example, Lemma 2.2.44 in [13]). 

 \medskip

 In Sections 3-5 we will use some elementary facts from the theory of orderable groups. We would like to quote the following well known folklore result.
   
 \medskip

 \begin{lem}\label{lem:orderable} Let $\Gamma $ be a left orderable countable group with a left-invariant order $<$. Then $\Gamma $ admits a faithful representation into $\mathrm{Homeo}_{+}(\mathbb{R})$ without a global fixed point, moreover, if $g_1, g_2\in \Gamma $ with $g_1 < g_2$ then $g_1(0) < g_2(0)$. 
\end{lem}
 
 \medskip

 The embeddability claim is stated and proved, e.g., in [13]; see Theorem 2.2.19 there. The ``moreover" part of the claim also follows immediately from the construction in the same proof.  

  \bigskip
  
  \section{Proof of Theorem \ref{thm:easy}}

\begin{lem} \label{lem:maximalorbital} Let $\Gamma \leq \mathrm{Homeo}_{+}(I)$ with $h(\Gamma ) < \infty $, and $(I_1,f_1), (I_2, f_2)\subset (0,1)$ be signed orbitals of $\Gamma $ where $I_2$ is maximal. Then either $I_1\subseteq I_2$ or $I_1\cap I_2 = \emptyset $.
\end{lem}

\medskip

{\bf Proof.} Let $I_1 = (a_1,b_1), I_2 = (a_2,b_2)$. We may assume that $f_i(x) > x$ for all $x\in I_i, 1\leq i\leq 2$. If neither of the conditions $I_1\subseteq I_2$ or $I_1\cap I_2 = \emptyset $ hold then,  since the interval $I_2$ is maximal, without loss of generality, we may assume that $a_1 < a_2 < b_1 < b_2$. 

\medskip

 If $b_2\notin Fix(f_1)$ then we consider the following two cases separately.
 
 \medskip

 {\em Case 1.} $f_1(b_2) < b_2$:  
 
  In this case, for every positive integer $n$, the element $f_1^nf_2f_1^{-n}$ has an orbital $(c_n,d_n)$ where $c_n = f_1^n(a_2), d_n = f_1^n(b_2)$, moreover, we have $c_1 < c_2 < c_3 < \ldots $ and $d_1 > d_2 > d_3 > \ldots $. This contradicts the assumption that the strict height of $\Gamma $ is finite.   
  
  \medskip
  
  {\em Case 2.} $f_1(b_2) > b_2$: 
  
  In this case, we shall consider the elements $g_n = f_1^{-n}f_2$, for $n\geq 1$. Let \begin{center} $p_n = \max \{x : a_2 < x < b_1, \ \mathrm{and} \ f_1^n(x) = f_2(x)\}$, \end{center} \ \begin{center} $z = \max \{x\in I_2 : f_1(x) = x\}$,\end{center} \ \begin{center} $q_n = \min \{x : z < x < b_2, \ \mathrm{and} \ f_1^n(x) = f_2(x)\}$. \end{center} 
        
 \medskip
       
  Notice that, by continuity, $p_n, z, q_n$ exist. Then for every positive integer $n$, $(p_n,q_n)$ is an orbital of the element $g_n$. Moreover, we have $p_1 < p_2 < p_3 < \ldots $ and $q_1 > q_2 > q_3 > \ldots $. Thus we obtain that $h_{strict}(\Gamma ) = \infty $ which contradicts the assumption.  
  
  \medskip
  
   Thus we established that $b_2\in Fix(f_1)$. Then $f_1^{-n}f_2f_1^n$ will have an orbital $(f_1^{-n}(a_2), b_2)$, and since $a_2\notin Fix(f_1)$ we obtain that $h(\Gamma ) = \infty $. $\square $
 
\bigskip

 \begin{prop} \label{prop:specialfixed} Every special subgroup $\Gamma \leq \mathrm{Homeo}_{+}(I)$ is  metabelian. Moreover, if in addition $h(\Gamma ) < \infty $, then $\Gamma $ is Abelian.
 \end{prop}
 
 \medskip
 
 The proof will follow the idea (from [2]) of the proof of the fact that a subgroup  of $\mathrm{Homeo}_{+}(I)$ where every non-identity element has at most one fixed point is necessarily metaabelian (originally due to Solodov (unpublished), later proved also by Barbot [4], Kovacevic [12], and Farb-Franks [9]). The issue here is to replace the condition ``every element has at most one fixed point" with ``every element is special", i.e. we do allow more than one fixed points (even infinitely many fixed points) for the elements of $\Gamma $, but we demand that no element has an inner orbital. For such a group $\Gamma $ we can introduce the following natural bi-invariant order: for $f,g\in \Gamma $, let $s = \max \{z\in [0,1] \ | \ f(x) = g(x), \forall x\in [0,z]\}$.  We will say $f\prec g$ if there exists an $\epsilon > 0$ such that $f(x) < g(x)$ for all $x\in (s, s+\epsilon )$. Then $\prec $ defines a bi-invariant order in $\Gamma $.
 
 \medskip

 {\bf Proof.} If all finitely generated subgroups of a group are metaabelian then the group is metaabelian.
Hence we may assume that $\Gamma $ is finitely generated with a fixed finite
generating set. Let $f$ be the biggest generator of $\Gamma $. Then $$\Gamma _f =\{g\in \Gamma \ | \ g^n\prec f, \ \mathrm{for \ all} \ n\in \mathbb{Z}\}$$ is a normal subgroup and $\Gamma /\Gamma _f$ is Archimedean therefore Abelian.

\medskip

 Without loss of generality, we may assume that $\Gamma $ is irreducible, i.e. it has no global fixed point in $(0,1)$. Since $f$ is the biggest generator, we may also assume that $Fix(f)\cap (0,\epsilon ) = \emptyset $ for some $\epsilon > 0$.
 
 \medskip

  Now, let $h\in \Gamma _f$ such that $h$ has at least one fixed point (if such $h$ does not exist then $\Gamma _f$ is Abelian, therefore $\Gamma $ is metaabelian) and $h$ is positive. We may also assume that $\Gamma _f$ has no global fixed point. (if $\Gamma _f$ has a global fixed point then by specialness of $\Gamma $ and by H\"older's Theorem, we obtain immediately that $\Gamma _f$ is Abelian, hence $\Gamma $ is metaabelian.)

 \medskip
 
 Let $s = max \{z\in [0,1] : h(x) = x, \forall x\in [0,z]\}$. If $s>0$, then $h$ does not fix any point in $(s,1)$. Also, we can conjugate $h$ to $h_1$ such that $0 < s_1 < s$ where $s_1 = max \{z\in [0,1] : h_1(x) = x, \forall x\in [0,z]\}$, and $h_1$ does not fix any point in the interval $(s_1,1)$. Then for sufficiently big integer $n$, $h_1(u) = h^n(u)$ for some $u\in (s,1)$. Then $h_1^{-1}h^n$ is not special; contradiction.
 
 \medskip  
 
 Assume now $s=0$. If $f$ has no fixed point, then for sufficiently big $n$, $fh^{-n}$ will not be special, and we again obtain a  contradiction. So we may assume that $f$ has a fixed point. Let $a$ be the minimal fixed point of $h$ in $(0,1)$. By conjugating $f$ by the element of $\Gamma _f$ if necessary, we may assume that the minimal fixed point of $f$ in $(0,1)$ is bigger than $a$. Then, again, then for sufficiently big $n$, $fh^{-n}$ will not be special. $\square $
 
 \bigskip
 
 \begin{rem}\label{rem:solodov} Proposition \ref{prop:specialfixed} generalizes Solodov's Theorem to the class of special subgroups. In [2] we have proved another generalization of Solodov's Theorem for subgroups of higher regularity where every element has at most $N$ fixed points. The idea of the above proof seems to be useful in obtaining similar results for subgroups of higher regularity where the condition ``special" is replaced with more general naturally extended condition ``$N$-special".
 \end{rem} 
 
 \medskip
 
 Proposition \ref{prop:specialfixed} immediately implies the following claim.
 
 \medskip

 \begin{lem}\label{lem:abelian} Let $\Gamma \leq \mathrm{Homeo}_{+}(I)$ be a non-Abelian subgroup with $h(\Gamma ) < \infty $. Then $\Gamma $ possesses a maximal orbital. $\square $ 
 \end{lem}
 
  \medskip
 
 {\bf Proof of Theorem \ref{thm:easy}.} We may assume that $\Gamma $ has no global fixed point. 
 
 \medskip
 
  The proof will be by induction on $N = h(\Gamma )$. For the base of induction, let $N = 1$. Then since $\Gamma $ has no global fixed point, we obtain that no non-identity element of $\Gamma $ has a fixed point. Then by H\" older's Theorem, $\Gamma $ is Abelian.
  
 \medskip  
  
   Now, assume that the claim holds for all subgroups with height less than $N\geq 2$, and assume that $h(\Gamma ) = N$. If $\Gamma $ is not Abelian then by Lemma \ref{lem:abelian} some element of $\Gamma $ contains a maximal orbital $I$. Then, by Lemma \ref{lem:maximalorbital}, for every $f\in \Gamma $ either $f(I)=I$ or $f(I)\cap I = \emptyset $. Let $$\Omega = \{J : J = f(I) \ \mathrm{for \ some} \ f\in \Gamma \},   \Gamma _0 = \{f\in \Gamma \ | \ f(J)=J \ \mathrm{for \ all} \ J\in \Omega\}.$$ Notice that $\Gamma _0$ is a normal subgroup with $h(\Gamma _0) = h(\Gamma )-1$. Thus by inductive hypothesis $\Gamma _0$ is solvable with solvability degree at most $N-1$. Also, for every $f\in \Gamma \backslash \Gamma _0$, the points 0 and 1 are the accumulation points of the set $\displaystyle \mathop{\sqcup }_{n\in \mathbb{Z}}f^{n}(I)$. Then $\Gamma /\Gamma _0$ is  Archimedian, hence Abelian. Thus $\Gamma $ is solvable of solvability degree at most $N$. $\square $   
 
 \medskip
 
 We also would like to remark that one can be more precise about the degree of solvability but this requires additional set of arguments.
       
 \bigskip
 
 \section{Proof of Theorem \ref{thm:hard}}
 
  We need to prove an analogue of Lemma \ref{lem:maximalorbital}. For this purpose, we need a notion of quasi-orbital.
  
   \begin{defn} An open interval $(a,b)$ will be called a quasi-orbital if there exists an infinite countable tower $T$ such that 
 
 1. the union of orbitals of $T$ equals $(a,b)$
 
 2. there exists a point $p\in \{a,b\}$ such that all orbitals of $T$ share the end $p$.
 
 3. no orbital of $T$ equals $(a,b)$.
 \end{defn}
 
  The end $p$ of the quasi-orbital $(a,b)$ will be called the {\em heavy end}.  A quasi-orbital is called inner if it does not share an end with the interval $(0,1)$. An inner quasi-orbital is called maximal if it is not properly contained in another inner quasi-orbital. 
 
 \medskip

 \begin{lem} \label{lem:analog} Let $\Gamma \leq \mathrm{Homeo}_{+}(I)$ with $h_{strict}(\Gamma ) < \infty $, and let $I_1, I_2\subset (0,1)$ be quasi-orbitals of $\Gamma $ where $I_2$ is maximal. Then either $I_1\subseteq I_2$ or $I_1\cap I_2 = \emptyset $.
\end{lem}

\medskip

{\bf Proof.} Let $I_1 = (c_1,b_1), I_2 = (c_2,b_2)$. We may assume that the ends $b_1, b_2$ are the heavy ends of the quasi-orbitals $I_1, I_2$ respectively (all other cases are similar). If neither of the conditions $I_1\subseteq I_2$ or $I_1\cap I_2 = \emptyset $ hold then, without loss of generality, we may assume that $c_1 < c_2 < b_1 < b_2$. Then we may assume that there exist two elements $f_1, f_2\in \Gamma $ with orbitals $(a_1, b_1), (a_2, b_2)$ such that $c_1 < a_1 < c_2 < a_2 < b_1 < b_2$, and $f_i(x) > x$, for all $x\in (a_i, b_i), 1\leq i\leq 2$.  

\medskip

 If $b_2\notin Fix(f_1)$ then we consider the following two cases separately.
 
 \medskip

 {\em Case 1.} $f_1(b_2) < b_2$:  
 
  In this case, for every positive integer $n$, the element $f_1^nf_2f_1^{-n}$ has an orbital $(t_n,d_n)$ where $t_n = f_1^n(a_2), d_n = f_1^n(b_2)$, moreover, $t_1 < t_2 < t_3 < \ldots $ and $d_1 > d_2 > d_3 > \ldots $. This contradicts the assumption that the strict height of $\Gamma $ is finite.   
  
  \medskip
  
  {\em Case 2.} $f_1(b_2) > b_2$: 
  
  In this case, we shall consider the elements $g_n = f_1^{-n}f_2, n\geq 1$. Let \begin{center} $p_n = \max \{x : a_2 < x < b_1, \ \mathrm{and} \ f_1^n(x) = f_2(x)\}$,\end{center} \ \begin{center} $z = \max \{x\in I_2 : f_1(x) = x\}$, \end{center} \ \begin{center} $q_n = \min \{x : z < x < b_2, \ \mathrm{and} \ f_1^n(x) = f_2(x)\}$.\end{center}

  \medskip
  
  Notice that, by continuity, $p_n, z, q_n$ exist. Then for every positive integer $n$, $(p_n,q_n)$ is an orbital of the element $g_n$, moreover, we have $p_1 < p_2 < p_3 < \ldots $ and $q_1 > q_2 > q_3 > \ldots $. This implies that $h_{strict}(\Gamma ) = \infty $ which contradicts the assumption.  
  
  \medskip
  
   Thus we established that $b_2\in Fix(f_1)$. 
  
   \medskip
 
 Because of the condition $f_1(b_2) = b_2$, the interval $(b_1, b_2)$ contains at least one orbital of $f_1$. If $f_2^{-1}(x) \leq f_1^{n}(x)$ for all $x\in (b_1, b_2)$ and for all positive $n$, then there exists $\epsilon > 0$ and $n\geq 1$ such that  $f_2f_1^n(x) > x$ for all $x\in (c_2-\epsilon , b_2)$. Since $b_2\in Fix(f_1)$, this contradicts maximality of $I_2$. 

 \medskip
 
 Assume now $f_2^{-1}(x) > f_1^{m}(x)$ for some $m\geq 1$, and $x\in I$ where $I\subseteq (b_1, b_2)$ is an orbital of $f_1$. Then there exist distinct points $u, v\in I$ such that the graphs of $f_1^m$ and $f_2^{-1}$ cross at the points $u, v$ and in no other point inside $(u,v)$ (hence $(u,v)$ is an orbital of $f_2f_1^m$). Then there exists a sequence $(n_k)_{k\geq 1}$ of strictly increasing positive integers such that the element $f_2f_1^{mn_k}$ has an orbital $(u_k, v_k)$ where $(u,v) \subset (u_1, v_1) \subset (u_2, v_2) \subset \ldots $ and the inclusions are strict at both ends. Then $h_{strict}(\Gamma ) = \infty $. Contradiction. $\square $

 \medskip
 
 Next, we need to observe that Proposition \ref{prop:specialfixed} immediately implies the following analogue of Lemma \ref{lem:abelian}. 
 
 \medskip
 
 \begin{lem}\label{lem:metaabelian} Let $\Gamma \leq \mathrm{Homeo}_{+}(I)$ be a subgroup such that $\Gamma $ is non-metaabelian and $h_{strict}(\Gamma ) < \infty $. Then $\Gamma $ possesses a maximal quasi-orbital. 
 \end{lem}
 
 \medskip
 
 Now we are ready to prove Theorem \ref{thm:hard}. The proof will be very similar to the proof of Theorem \ref{thm:easy}; we will use Lemma \ref{lem:analog} and Lemma \ref{lem:metaabelian} (instead of Lemma \ref{lem:maximalorbital} and Lemma \ref{lem:abelian}). 
 
 \medskip
 
 {\bf Proof of Theorem \ref{thm:hard}.} We may assume that $\Gamma $ has no global fixed point. The proof is again by induction on $N = h_{strict}(\Gamma )$. For the base of induction, let $N = 1$. Then since $\Gamma $ has no global fixed point, we obtain that $\Gamma $ is special. Then by Proposition \ref{prop:specialfixed}, $\Gamma $ is metaabelian.
  
 \medskip  
  
   Now, assume that the claim holds for all subgroups with strict height less than $N\geq 2$, and assume that $h_{strict}(\Gamma ) = N$. If $\Gamma $ is not metaabelian then by Lemma \ref{lem:metaabelian} some element of $\Gamma $ contains a maximal quasi-orbital $I$. Then, by Lemma \ref{lem:analog}, for every $f\in \Gamma $ either $f(I)=I$ or $f(I)\cap I = \emptyset $. Let $$\Omega = \{J : J = f(I) \ \mathrm{for \ some} \ f\in \Gamma \},   \Gamma _0 = \{f\in \Gamma \ | \ f(J)=J \ \mathrm{for \ all} \ J\in \Omega\}.$$ Notice that $\Gamma _0$ is a normal subgroup with $h_{strict}(\Gamma _0) = h_{strict}(\Gamma )-1$. Thus by inductive hypothesis $\Gamma _0$ is solvable with solvability degree at most $N$. Also, for every $f\in \Gamma \backslash \Gamma _0$, the points 0 and 1 are the accumulation points of the set $\displaystyle \mathop{\sqcup }_{n\in \mathbb{Z}}f^{n}(I)$. Then $\Gamma /\Gamma _0$ is  Archimedian, hence Abelian. Thus $\Gamma $ is solvable of solvability degree at most $N+1$. $\square $

 \vspace{1cm}
 
 \section{A solvable group with a bi-infinite strict tower} 
 
 Theorem \ref{thm:hard} naturally leads to a converse question of whether or not a subgroup with finite strict height is necessarily solvable. Let us recall that it is quite easy to find solvable subgroups with infinite height, see Remark \ref{rem:hvsh'}.  
 
 \medskip
 
 In this section we present an example of a 2-generated solvable group with infinite strict height thus proving Theorem \ref{thm:example}   
 
 \medskip
 
 Let $\Gamma $ be a group generated by three elements $t, a, b\in \Gamma $. Let us assume that the following conditions hold:
  
  \medskip
  
  (i) $\Gamma $ is solvable,
  
   \medskip
   
  (ii) $\Gamma $ is left-orderable with a left order $<$, moreover, $$b^{-1} < a^{-1} < t^{-1} < 1 < t < a < b,$$
  
   \medskip
   
   (iii) $t^mat^{-m} < t^nat^{-n}$ for all integers $m < n$.
     
   \medskip
   
  To state the last two conditions we need to introduce some notations: let $C$ denotes the cyclic subgroup of $\Gamma $ generated by $t$, $G$ denotes the subgroup generated by $t$ and $a$.
    
  \medskip
  
  (iv) if $g\in C, f\in \Gamma \backslash C, 1 < f$ then $f^{-1} < g^{-1} < 1 < g < f$;
    
  \medskip
  
  (v) if $g\in G,  f\in \Gamma \backslash G, 1 < f$ then $f^{-1}< g^{-1} < 1 < g < f$ .
   
   \medskip

   We are postponing the construction of $\Gamma $ with properties (i)-(v) till the end. 
    
   \bigskip
   
   Let us now observe some implications of conditions (i)-(v):
   
   \medskip
   
   Because of (ii), by Lemma \ref{lem:orderable}, $\Gamma $ is embeddable in Homeo$_{+}(\mathbb{R})$. Moreover, we can embed $\Gamma $ faithfully in Homeo$_{+}(\mathbb{R})$ such that the following conditions hold:

   \medskip
   
   (c1) if $g_1, g_2\in \Gamma , g_1 < g_2$ then $g_1(0) < g_2(0)$ (in particular, $g(0) > 0$ for all positive $g\in \Gamma )$;
   
   \medskip
   
   (c2) $\Gamma $ has no fixed point.
  
  \medskip
  
  We intend to show that if all the conditions (i)-(v), (c1)-(c2) are satisfied, then the subgroup $\Gamma $ contains a strict infinite tower.

  \medskip
  
  For any $g\neq 1$, by condition (c1), the set $Fix(g) \cap (0,\infty )$ is either empty or contains a minimal element; in the latter case, let $F_{+}(g)$ denotes this minimal fixed point. Similarly, the set $Fix(g) \cap (-\infty ,0)$ is either empty or contains a maximal element; in the latter case, let $F_{-}(g)$ denotes this maximal fixed point.

  \medskip
  
  Notice that $a^{-1}(0) < 0 < a(0)$. If $Fix(a)\cap (0,\infty ) = \emptyset $ then $a^n(0) > b(0)$ for a sufficiently big $n$ which contradicts the conditions (ii), (v) and (c1). Thus $a$ has a fixed point in the interval $(0, \infty )$. Similarly (by comparing $a^{-n}$ with $b^{-1}$) we obtain that $a$ has a fixed point in $(-\infty , 0)$ as well.   
  
  \medskip
  
  Let $p = F_{+}(a), q = F_{-}(a)$. If $p\in Fix(t)$, we have $t^{-n}at^n(0) > a(0)$ for a sufficiently big positive $n$ but this contradicts condition (iii). Hence, $p\notin Fix(t)$. Moreover, if $t(p) > p$ then, again, $t^{-1}a^nt(0) > a(0)$ for sufficiently big $n$ which contradicts condition (iii). Thus $t(p) < p$.  Then we have $$0 < F_{+}(t^{-m}at^m) < F_{+}(t^{-n}at^n) < p \ \mathrm{for \ all \ positive} \ m > n.$$ Similarly, we obtain  $$q < F_{-}(t^{-n}at^n) < F_{-}(t^{-m}at^m) < 0 \ \mathrm{for \ all \ positive} \ m > n.$$ Thus we obtain a bi-infinite strict tower $\{((F_{-}(t^{-k}at^k), F_{+}(t^{-k}at^k)), t^{-k}at^k)\}_{k\in \mathbb{Z}}$.

  \bigskip

  {\bf Construction of} $\Gamma $: Let us now construct $\Gamma $ with properties (i)-(v).

 \medskip
   
  We consider the rings $T = A = \mathbb{Z}, B = \mathbb{Z}[\frac{1}{2}]$. We will identify $t, a, b$ with the identity elements of the rings $T, A, B$ respectively. 
  
  \medskip
  
  We let $\Gamma  = B\wr (A\wr T)$. Let also $D = \oplus _{i\in \mathbb{Z}}A_i$ where $A_i, i\in \mathbb{Z}$ is an isomorphic copy of the ring $A$. Similarly, let $\Omega  = \oplus _{i\in \mathbb{Z}}H_i$, where $H_i, i\in \mathbb{Z}$ is an isomorphic copy of the ring $B$. An element of $\Omega $ can be represented by a vector ${\bf x} = (\ldots , x_{-1}, x_0, x_1, \ldots )$ where all but finitely many coordinates are zero. 
  
  \medskip

  The group $A\wr T = \mathbb{Z}\wr \mathbb{Z}$ acts on $\Omega $ as follows:  
  
  for all ${\bf x} = (\ldots , x_{-1}, x_0, x_1, \ldots )\in \Omega $, 
  
  $t({\bf x}) = {\bf y}$ where ${\bf y} = (\ldots , y_{-1}, y_0, y_1, \ldots ), y_{n} =x_{n-1}, \forall n\in \mathbb{Z}$. (so $t$ acts by a shift);
  
  $a({\bf x}) = {\bf y}$ where ${\bf y} = (\ldots , y_{-1}, y_0, y_1, \ldots ), y_{n} = 2x_n, \forall n\in \mathbb{Z}$. 
  
  \medskip
  
  Then $\Gamma = (A\wr T)\ltimes \Omega $ is the semidirect product with respect to the described action. 
   
   \medskip
   
   By construction, $\Gamma $ is solvable. To discuss conditions (ii)-(v), let us recall a basic fact about left-orderable groups.
   
   \medskip
   
  \begin{lem} Let a group $G_1$ acts on a group $G_2$ by automorphisms. Let $\prec _1, \prec _2$ be left orders on $G_1, G_2$ respectively, and assume that the action of $G_1$ on $G_2$ preserves the left order (i.e. if $g\in G_1, x_1, x_2\in G_2, x_1 \prec _2 x_2$ then $g(x_1)\prec _2 g(x_2)$).
   
   Then there exists a left order $<$ in $G_1\ltimes G_2$ which satisfies the following conditions:

   \medskip
   
   {\bf 1)} if $g_1, f_1\in G_1, g_1\prec _1 f_1$ then $(g_1,1) < (f_1, 1)$;
   
   \medskip
   
   {\bf 2)} if $g_2, f_2\in G_2, g_2\prec _2 f_2$ then $(1,g_2) < (1, f_2)$; 
   
   \medskip
   
   {\bf 3)} if $g_1\in G_1\backslash \{1\}, g_2\in G_2\backslash \{1\}, 1\prec _2 g_2$, then $(g_1, 1) < (1, g_2)$.
   \end{lem}
   
 \medskip
 
  {\bf Proof.} We define the left order on the semidirect product $G_1\ltimes G_2$ as follows: given $(g_1, f_1), (g_2, f_2)\in G_1\ltimes G_2$ we define $(g_1, f_1) < (g_2, f_2)$ if and only if either $f_1 \prec _2 f_2$ or $f_1 = f_2, g_1\prec _1 g_2$. Then the claim is a direct check. $\square $
  
  \bigskip
  
  The left order $<$ on the semidirect product $G_1\ltimes G_2$ constructed in the proof of the lemma will be called the {\em extension of} $\prec _1$ {\em and} $\prec_2$.
  
  \bigskip
  
  Let us now introduce a left order $\prec _1$ on the additive subgroup of the ring $B$. Notice that the additive groups $A, B$ are subgroups of $\mathbb{R}$, and we define $\prec _1$ on $B$ to be simply the restriction of the natural order of $\mathbb{R}$.

  \medskip
  
  Now we introduce an order $\prec _2$ on $D$. Let ${\bf x} = (\ldots , x_{-1}, x_0, x_1, \ldots ), {\bf y} = (\ldots , y_{-1}, y_0, y_1, \ldots )\in D$. We say ${\bf x}\prec _2 {\bf y}$ if and only if  $$\mathrm{min}\{k \ | \ x_k < y_k\} < \mathrm{min}\{k \ | \ y_k < x_k\}.$$ Similarly, we introduce an order $\prec _3$ on $\Omega $. Let ${\bf x} = (\ldots , x_{-1}, x_0, x_1, \ldots ), {\bf y} = (\ldots , y_{-1}, y_0, y_1, \ldots )\in \Omega $. We say ${\bf x}\prec _3 {\bf y}$ if and only if  $$\mathrm{min}\{k \ | \ x_k < y_k\} < \mathrm{min}\{k \ | \ y_k < x_k\}.$$
  
  \medskip
  
  Notice that $A\wr T$ is isomorphic to $T\ltimes D$ (where the semidirect product is with respect to the standard action of $T$ on $D$, by a shift, i.e. for all vectors ${\bf x} = (\ldots , x_{-1}, x_0, x_1, \ldots )\in D$, we have $$t({\bf x}) = {\bf y}$$ where ${\bf y} = (\ldots , y_{-1}, y_0, y_1, \ldots ), y_{n} =x_{n-1}, \forall n\in \mathbb{Z}$). Notice that the action of $T$ on $D$ preserves the left order $\prec _2$. Then, let $\prec _4$ be the extension of $\prec _1$ and $\prec _2$. Having the left order $\prec _4$ on $A\wr T$, we define the left order $<$ on $\Gamma = BS(1,2)\ltimes \Omega $ to be the extension of $\prec _4$ and $\prec _3$ [again, notice that the action of $A\wr T$ on $\Omega $ preserves the left order $\prec _3$]. The group $\Gamma = (A\wr T)\ltimes \Omega $ with the left order $<$ satisfies conditions (ii)-(v). 

\medskip

 It remains to notice that the subgroup $\Gamma _1$ generated by $t$ and $a$ is metaabelian (indeed it is isomorphic to $\mathbb{Z}\wr \mathbb{Z}$). $\square $     
  
  \medskip
  
  \begin{rem} Notice that the elements $t^nat^{-n}, n\in \Z$ all commute thus we have an Abelian subgroup $\Gamma _0\leq \Gamma $ with infinite strict height. This Abelian subgroup is not finitely generated. In fact, it is not difficult to show that the strict height of a finitely generated Abelian subgroup of $\mathrm{Homeo}_{+}(I)$ is always finite.
We also would like to point out that $\Gamma _0$ lies totally in a subgroup generated by $t$ and $a$ only; the role of adding an extra generator $b$ in the construction is to find a left-invariant order which helps to embed the group $\Gamma $ in a special way described in the construction. With a somewhat different construction, one can find a more direct embedding of the group $\Z \wr \Z$ into $\mathrm{Homeo}_{+}(I)$, with an Abelian subgroup of inifnite strict height. 
  \end{rem}
  
  \medskip
  
  \begin{rem} The construction of $\Gamma $ has no realization in $C^2$ regularity. This is a direct consequence of Koppel's Lemma [11] and the fact that the elements $t^nat^{-n}, n\in \Z$ all commute forming a bi-infinite strict tower.
  \end{rem}
  
 \bigskip
 
 \section{Proof of Theorem \ref{thm:analytic}}
 
  In this section we will study the height and the strict height of subgroups of analytic diffeomorphisms of $I$.
  
  \medskip
  
  Notice that by Remark \ref{rem:hvsh'}, a solvable (even a metaabelian) subgroup of $\mathrm{Diff}^{\omega }_{+}(I)$ may have infinite height. What about strict heights? We prove the following proposition which is a reformulation of Theorem \ref{thm:analytic}.  
 
 \begin{prop}\label{prop:analytic} A solvable subgroup of Diff $_{\omega }^{+}(I)$ does not possess a strict tower of length two.
 \end{prop}
 
 \medskip
 
 In the proof we will use the following lemma which is interesting in itself. 
 
 \medskip
 
 \begin{lem}\label{lem:discrete} Let $\Gamma \leq \mathrm{Diff}^{\omega }_{+}(I)$ be a non-Abelian subgroup with a non-trivial Abelian normal subgroup. Then $\Gamma $ is not $C_0$-discrete. Moreover, any non-trivial Abelian normal subgroup of $\Gamma $ is not $C_0$-discrete. 
 \end{lem}   
 
 \medskip
 
 {\bf Proof.} We may assume that $\Gamma $ is irreducible. Let $N$ be a non-trivial Abelian normal subgroup of $\Gamma $. Then we claim that $N$ acts freely. Indeed, if a non-identity element $f$ of $N$ has a fixed point $p\in (0,1)$, then any other element $f_1$ fixes $p$ (because otherwise $f$ has infinitely many fixed points). Then $N$ is not irreducible. Then, by irreducibility of $\Gamma $, any element of $N$ has infinitely many fixed points; contradiction.
 
 \medskip
 
 Thus we established that the action of the subgroup $N$ on $(0,1)$ is free. Let $\epsilon > 0$ and $h\in N\backslash \{1\}$ where $h(x)>x$ for all $x\in (0,1)$. Let also  $g\in \Gamma $ such that $g$ does not act freely and $g(x) < x$ near 1. Then $[g,h]\neq 1$. By replacing $h$ with $[g,h]$ if necessary, we may assume that $h'(0) = h'(1) = 1$. 
 
 \medskip
 
 We have either $ghg^{-1}(x) < h(x), \forall x\in (0,1)$ or $g^{-1}hg(x) < h(x), \forall x\in (0,1)$. By replacing $g$ with its conjugate if necessary, we may also assume that $\max (Fix(g)\cap (0,1)) \in (0,\epsilon )$ and $g(x) < x$ on $(\epsilon , 1)$.  Then $ghg^{-1}(x) < h(x), \forall x\in (0,1)$; and $x < g^nhg^{-n}(x) < g^mhg^{-m}(x)$ for all $x\in (0,1)$ and for all $0 < m < n$. 
 
 \medskip
 
 Notice that $\Gamma $ always contains a crossed pair, hence a free semigroup (since $h$ acts freely, for sufficiently big $n$, either the pair $(h^ngh^{-n}, h^{-n}gh^n)$ or the pair $(h^ngh^{-n}, h^{-n}g^{-1}h^n)$ is crossed). Therefore, if $g'(0) = g'(1) = 1$ then by Theorem A in [3] (more precisely, by the proof of Theorem A), $[\Gamma , \Gamma ]$ is not $C_0$-discrete. \footnote{in Theorem A, as stated in [3], the claim is that if $[\Gamma ,\Gamma ]$ contains a free semigroup then $\Gamma $ is not discrete, but notice that, first, the only property of the commutator subgroup used in the proof is the fact that an element of $[\Gamma , \Gamma ]$ has derivative 1 at the endpoints of the interval $[0,1]$; second, the elements constructed in the proof which are arbitrarily close to identity in $C_0$-metric are indeed from the subgroup $[\Gamma , \Gamma ]$. Thus in [3], we have indeed proved the following claim: {\em if $\Gamma $ contains a free semigroup on generators $f_1, f_2$ and $f_i'(0) = f_i'(1) = 1, i = 1, 2$ then the commutator subgroup $[\Gamma , \Gamma ]$ is not $C_0$-discrete.}} 
 
 So we may assume that $g'(1) \neq 1$. Let also $p = \max (Fix(g)\cap (0,1))$. 
 Then either $g^nhg^{-n}(x)\to x$ as $n\to \infty $ or, more generally, there exists a non-decreasing function $\phi :[0,1]\to [0,1]$ such that $\phi (0) = 0, \phi (1) = 1$ and $g^nhg^{-n}(x)\to \phi (x)$ as $n\to \infty $. However, we can make a stronger claim about the function $\phi $: first, notice that $\phi (p) = p$. On the other hand, by $C^2$-differentiability, for all $x\in (p,1)$ there exists a $\delta > 0$ such that the sequence $(g^nhg^{-n})'(z)$ is bounded on the interval $(x-\delta , x+\delta )$, i.e. there exists an $M > 1$ such that $\frac{1}{M} < (g^nhg^{-n})'(z) < M$ for all $z\in (x-\delta , x+\delta )$. To see this, notice that for sufficiently small $\delta > 0$, the monotone sequence $a_n = g^{-n}(z), n\geq 1$ tends to 1 as $n\to \infty $. Then, since $g'(1) \neq 1, h'(1) = 1$, for sufficiently large $n$, we have $h(a_n)\in (a_n, a_{n+1})$. On the other hand, the quantity $(a_n-a_{n+1})$ tends to zero exponentially hence faster than $\frac{1}{n}$. Then by $C^2$-differentiability of $g$, for sufficiently big $n$ we have $$\displaystyle \max _{u,v\in [a_n,a_{n+1}]}\frac{g'(u)}{g'(v)} \leq 1+\frac{1}{n}.$$ Since the product $(1+\frac{1}{n})^n$ is bounded, the claim about the boundedness of the derivative $(g^nhg^{-n})'(z)$ follows immediately from the chain rule. Then the function $\phi $ is a homeomorphism on the interval $[p,1]$; hence, for all $x\in [p,1]$, $(g^nh^{-1}g^{-n})g^mhg^{-m}(x)\to x$ as $m, n\to \infty $. Since $\epsilon $ is arbitrary, we obtain that $\Gamma $ is not $C_0$-discrete. It remains to observe that $g^nhg^{-n}\in N$ for every $n\in \mathbb{Z}$ and for every $g\in \Gamma $. Thus $N$ is not $C_0$-discrete. $\square $
 
 \medskip
 
  \begin{rem}\label{rem:discrete} As a corollary of Lemma \ref{lem:discrete} we obtain that a solvable non-Abelian subgroup of $\mathrm{Diff}^{\omega }_{+}(I)$ is never $C_0$-discrete. On the other hand, it is proved in [3] that a non-solvable (non-metaabelian) subgroup of $\mathrm{Diff}^{1+\epsilon }_{+}(I)$ (of $\mathrm{Diff}^{2}_{+}(I)$) is never $C_0$-discrete. Thus, in the context of analytic diffeomorphisms, Lemma \ref{lem:discrete} gives us a new result not covered by the results of [3]. 
  \end{rem}     
 
  \medskip
  
 {\bf Proof of Proposition \ref{prop:analytic}.} The claim of the proposition follows from the fact that for any irreducible solvable subgroup $\Gamma $, every non-identity element has at most one fixed point in $(0,1)$. Indeed, let $\Gamma $ be a solvable subgroup of Diff $_{\omega }^{+}(I)$. By the result of E.Ghys [10], $\Gamma $ is metaabelian. Then any non-identity element of $[\Gamma , \Gamma ]$ acts freely. Then by Lemma \ref{lem:discrete}, the group $[\Gamma ,\Gamma ]$ is either trivial or not $C_0$-discrete. 
 
 \medskip
 
 If $[\Gamma ,\Gamma ]$ is trivial then $\Gamma $ is Abelian hence any point which is fixed by a non-identity element is fixed by the whole group $\Gamma $. By irreducibility, $\Gamma $ acts freely, hence $h_{strict}(\Gamma ) = 1$. 
 
 \medskip
 
 If $[\Gamma , \Gamma ]$ is not discrete then let $f\in \Gamma \backslash \{1\}$ with some two consecutive fixed points $a, b\in (0,1)$ where $a<b$. By non-discreteness of $[\Gamma , \Gamma ]$, there exists an element $h\in [\Gamma , \Gamma ]\backslash \{1\}$, such that $a \leq h(a) < b$. On the other hand, by freeness of the action of $[\Gamma , \Gamma ]$, $h(a)\neq a$. Thus $hfh^{-1}$ has a fixed point in the open interval $(a,b)$. Then the commutator $[h,f]$ is non-trivial and has a fixed point. Then the action of $[\Gamma , \Gamma ]$ is not free. Contradiction. $\square $ 

   \vspace{1cm}

 {\em Acknowledgment:} We would like to thank the anonymous referee for useful remarks and suggestions which helped us improve the paper.
  
  \vspace{3cm}
  
  {\bf R e f e r e n c e s:}
  
  \bigskip
  
  [1] Akhmedov, A. \ Girth alternative for subgroups of $PL_o(I)$. \ Preprint. http://arxiv.org/pdf/1105.4908.pdf

\medskip

  [2] Akhmedov, A. \ Extension of H\"older's Theorem in $\mathrm{Diff}_{+}^{1+\epsilon }(I)$. \ Preprint. \  http://arxiv.org/abs/1308.0250.pdf 
  
  \medskip
  
  [3] Akhmedov, A. \ A weak Zassenhaus Lemma for discrete subgroups of Diff(I). {\em Algebraic and Geometric Topology,} {\bf 14}, no. 1, (2014) 539-550, \ http://arxiv.org/pdf/1211.1086.pdf

\medskip

  [4] T.Barbot, \ Characterization des flots d'Anosov en dimension 3 par leurs feuilletages faibles, \ {\em Ergodic Theory and Dynamical Systems} {\bf 15} (1995), no.2, 247-270.
  
  \medskip
  
  [5] Bleak, C. \ A geometric classification of some solvable groups of homeomorphisms. {\em Journal of London Mathematical Society} (2), {\bf 78}, (2008) no. 2, 352-372. 

\medskip
  
  [6] Bleak, C. \ An algebraic classification of some solvable groups of homeomorphisms. \ {\em Journal of Algebra.} {\bf 319}, 4, p. 1368-1397.

\medskip
  
  [7] Bleak, C. \ Solvability in Groups of Piecewise Linear Homeomorphisms of the Unit Interval. \ Ph.D.  Thesis. \ SUNY Binghamton, 2005.

\medskip

 [8] Burslem, L, Wilkinson, A. \ Global rigidity of solvable group actions on $\mathbb{S}^1$. \ {\em Geometry and Topology}, {\bf 8},  (2004), 877-924.
 
 \medskip
 
 [9] Farb, B., Franks, J. \ Groups of homeomorphisms of one-manifolds, II: Extension of Holder's Theorem. {\em Transactions of the  American Mathematical Society.} {\bf 355}, (2003) , no.11, 4385-4396. 
  
  \medskip

  [10] Ghys, E. \ Sur les groupes engendrés par des difféomorphismes proches de l'identité.  Bol. Soc. Brasil. Mat. (N.S.) 24 (1993), no. 2, 137-178.
  
  \medskip
  
  [11] Kopell, N. Commuting Diffeomorphisms. \ In Global Analysis, Pros. Sympos. Pure Math XIV Berkeley, California, (1968) 165-184.
	  
	 \medskip
	 
	[12] Kovacevic, N. \ M\"obius-like groups of homeomorphisms of the circle. \ {\em Trans. Amer. Math. Soc.} {\bf 351} (1999), no.12, 4791-4822.
	
	\medskip
	
	[13] Navas, A. \ Groups of Circle Diffemorphisms. \ Chicago Lectures in Mathematics, 2011. {\em http://arxiv.org/pdf/0607481}

\end{document}